\documentclass[11pt]{article}
\usepackage{amsbsy,amsmath,amsthm,amsfonts,amssymb,amstext,amscd,longtable,flafter,float}

\newtheorem{theo}{\bf Theorem}[section]

\newtheorem{defi}{\bf Definition}[section]
\newtheorem{nota}{\bf Notation}[section]
\newtheorem{lem}{\bf Lemma}[section]

\begin{document}

\begin{center}
{\bf {\Large  Aspects of optimality of plans orthogonal through other factors}}

\vskip5pt

{\bf {\large Bhaskar Bagchi and Sunanda Bagchi*
\\[0pt]
1363, 10th cross, Kengeri Satellite Town, \\
Bangalore 560060, \\
India\\}}
\end{center}

\vskip10pt

\vskip5pt {\bf {\large Abstract  }}


\vskip5pt

 The concept of orthogonality through the block factor (OTB), defined  in Bagchi (2010), is extended here to
 orthogonality through a set (say S) of other factors. We discuss the impact of such an orthogonality
on the precision of the estimates as well as on the inference procedure. Concentrating on the case when
$S$ is of size two, we construct a series of plans in each of which every pair of other factors is orthogonal through
a given pair of factors.


Next we concentrate on plans through the block factors (POTB).  We construct POTBs for symmetrical experiments with two and three-level factors. The
plans for  two  factors are E-optimal, while those for three-level  factors are  universally optimal.
 Finally, we construct  POTBs for  $s^t(s+1)$ experiments, where $s \equiv 3 \pmod 4$ is a prime power. The plan is universally optimal.

\vskip10pt




\section{Introduction}
Morgan and Uddin (1996) have pioneered the path of deviation from the traditional condition of orthogonality among all factors, constructing plans where treatment factors are orthogonal to each other, but not necessarily orthogonal to the nuisance factors.
Later Mukherjee, Dey and Chatterjee (2002) discussed and constructed main effect plans (MEPs) on small-sized blocks with treatment factors non-orthogonal to the block factor, but orthogonal among themselves. The  plans constructed in both of these papers  satisfy optimality.
Das and Dey (2004) constructed plans with similar property. Restricting to blocks of size two, Bose and Bagchi (2007)
provided  plans satisfying  properties similar to those of the plans of Mukherjee, Dey and Chatterjee (2002), but with fewer blocks.
The condition of `orthogonality through the block factor' between a pair of  treatment factors was formally defined in Bagchi (2010). A plan in which every  treatment factors is orthogonal to every other one through the block factor is
 named  `plans orthogonal through the block factor' (POTB). By now,  new classes of POTBs have been constructed by many authors.
 Jacroux and his co-authors (20011, ... 2017) have come up with a number of such plans, mostly for two-level factors, many of them satisfying
optimality properties. Other authors include Chen,  Lin,  Yang,  and Wang, H.X (2015) and Saharay and Dutta 2016).

        In this paper we  define the concept of `orthogonality through a set of factors'. We show how such a property leads to the achievement of simplicity and precision in the data analysis. Next we go for  construction. We construct the following.
         (a) a series plans orthogonal through a pair of factors, (b) a few series of  plans orthogonal through the block factor, satisfying optimality property.

    In Section 2 we present definitions and notations for the set up. In Section 3 we define and discuss `orthogonality through a set of factors'. We have shown that regarding inference on a factor, say $A$, one may forget all factors other than those in a set $T$ if and only if $A$ is orthogonal to all factors through $T$ [see Theorem \ref {EstSSQ}]. In Section 4 we concentrate on construction. In Section 4.1 we construct an infinite series of plans orthogonal through a pair  of factors  [see Theorem \ref {2orthPlan}. In Section 4.2 we construct POTBs. Specifically, we obtain (a) an infinite series of E-optimal POTBs for two-level factors [see Theorem \ref {2power7n}], (b) an infinite series of universally optimal POTBs for three-level factors [see Theorem \ref {3^{3N}}] and (c)  an infinite series of universally optimal POTBs, for an asymmetrical experiment with bigger sets of levels [see Theorem \ref {level(s+1)}].

\section{Preliminaries}

We shall consider main effect plans (say ${\cal P}$) for an
experiment in which a block factor may or may not be present.

\begin{nota} \label{expt}
(a) The number of factors is denoted by $m$ and the number of runs
by $n$.

(b) ${\cal F}$ will denote the set of all factors of  ${\cal P}, S_A$ the set of levels of $A$ and $s_A = |S_A|, A \in {\cal F}$.
  We shall view the general effect as a factor, say $G$, so that $s_G = 1$.

(c) For  $x_A \in S_A, A \in {\cal F}$, the vector $x =(x_A : A \in {\cal F} )'$ represents a  level combination or run, in which $A$
is at level $x_A, A \in {\cal F}$.


(d) Fix $A : A \in {\cal F}$. The replication number $r^A (t)$ of
the level $t$ of factor $A$  is the number of runs in ${\cal P}$
in which $A$ is at level $t$. The replication vector $r_A$ is the $s_A \times 1$ vector
 with $r^A (t)$ as the $t$-th entry, $t \in S_A$. $ R_A$ denotes the diagonal matrix with  diagonal entries same
as those of $ r_A$ in the same order.

(e) For $A,B \in {\cal F}$, the $A$ versus $B$ incidence matrix is the $ s_A \times s_{B}$ matrix $N_{AB}$. The $(p,q)$th entry of
this matrix is  $n^{AB} (p, q)$, which is the number of runs $x \in
{\cal P}$ such that $x_A = p$ and $x_B = q$, where $p \in S_A, \; q \in s_B$. When $B = A, N_{AB} = R_A$.
\end{nota}

\begin{nota} \label{model}
(a) $ 1_n$ will denote the $n \times 1$ vector of all-ones, while $
J_{m \times n }$ will denote the $m \times n$ matrix of all-ones. We shall write $J_m$ for $J_{m \times m}$.

(b) For any $m \times n$ matrix $M$, ${\cal C} (M)$ will denote the
column space of $A$.

(c) The $n \times 1 $ vector of responses will be denoted by $Y$.

(d) The  $ s_A \times 1$ vector $\alpha^A$ will denote the vector of
unknown effects of $A, \;A \in {\cal F}$.
\end{nota}

 The model is expressed in matrix form as
\begin{equation}\label{modelEq} {\mathbf Y} =  {\mathbf X}  \alpha +\epsilon :
{\mathbf X} = [ {\mathbf X}_A : A \in {\cal F}], \; \;  \alpha = [\alpha^A,  A \in {\cal F}], \; \; \epsilon \sim N_n(0,\sigma^2 I_n).
\end{equation}
Here $Z \sim N_p(\mu,\Sigma)$ means that $Z$ is a random variable
following $p$-variate normal distribution with mean $\mu$ and
covariance matrix $\Sigma$. For $A \in {\cal F}, \; X_A$ is the
{\bf design matrix} for $A$. Thus,
${\mathbf X}_A $ is the $n \times s_A$ matrix having the $(u,t)$th
entry  1 if in the uth run the factor $A$ is set at level t and 0
otherwise, $1 \leq u \leq n,\; t \in s_A$. In particular, $X_G = 1_n$. The
following relations are well-known.

\begin{equation}\label{designMatrix}
{\cal C} (X_G) \subseteq {\cal C} (X_A), \; N_{AG} = r_A \mbox{ and } N_{AB} = X'_A X_B, \;
A,B \in {\cal F}.
\end{equation}

\vspace{.5em}

We shall use the following notations for the sake of compactness.
\begin{nota} \label{C-matandQ}
(a) For any $m \times n$ matrix $M$,  $P_M$ will denote the projection operator on the column space of $M$. Thus, $P_M =M (M'M)^-M'$, where $B^-$ denotes a g-inverse of $B$.

(b)  Let   $T$ be a  subset of ${\cal F}$.

 (i) $ X_T$ will denote $[  {\mathbf X}_A : A \in T]$.   Moreover, for an $A \in {\cal F}, \; X'_A X_T$
 will be denoted by $N_{AT}$ (which is consistent with (\ref  {designMatrix}).

  (ii) $\alpha^T$ will denote  $ [(\alpha^A)' : A \in T]'$. $\widehat{\alpha^T}$ will denote the least square
estimate of $\alpha^T$.

(iii)  $ P_A$ will denote the projection operator onto the column
space of $ X_A, A \in {\cal F}$. Further, $ P_T$ will denote the projection
operator onto the column space of $ X_T$.

(c) Consider a pair of disjoint subsets $T$ and $U$ of ${\cal F}$.  We
define the matrix ${\mathbf C}_{T;U}$ and the vector ${\mathbf Q}_{T;U}$ as follows.
\begin{eqnarray} \label{CQSS}
 {\mathbf C}_{T;U} & = & (( C_{AB;U}))_{A,B \in T},\;
  C_{AB;U}  =  {\mathbf X}'_A (I - P_U) {\mathbf X}_{B}, \\
{\mathbf Q}_{T;U} & = & (( Q_{A;U} ))_{A \in T}, \;  Q_{A;U} = {\mathbf X}'_A (I - P_U) {\mathbf Y}.\end{eqnarray}

(d) Let   $\bar{A} = {\cal F} \setminus \{A\}$. For the sake of simplicity we shall use the notation
${\mathbf C}_A$ (respectively ${\mathbf Q}_A$) instead of ${\mathbf C}_{A;\bar{A}}$ (respectively ${\mathbf Q}_{A;\bar{A}}$).

(e) Sum of squares : Fix a set of factors $T$. For $A  \notin T$, the sum of squares for $A$, adjusted for the factors
in $T$ will be denoted by  $SS_{A;T}$.

More generally  the combined sum of squares for the set of factors
$U$  adjusted for the set of factors $T$, (T disjoint from U) will be denoted by $SS_{U;T}$.
 \end{nota}

 \vspace{.5em}

{\bf Remark 3.1:} The matrix  ${\mathbf C}_A$ in (d) of Notation \ref {C-matandQ}
is referred to as the ``information matrix" or ``C-matrix" of $A$.
In order that every main effect contrast of $A$ is estimable, rank of ${\mathbf C}_A$ must be $s_A - 1$.
We, therefore, consider only the plans satisfying $Rank( {\mathbf C}_A) = s_A - 1$, for every $A \in {\cal F}$.
We shall refer to such a plan as `connected'.

\vspace{.5em}

The following relations are well-known.
\begin{lem} \label{NEandSS}
(a) The normal equation for the least square estimates of the vector of all effects is
 \begin{equation}  \label{NEalpha}
  {\mathbf X}'  {\mathbf X} \widehat{\alpha} = {\mathbf X}' {\mathbf Y}.
  \end{equation}

(b) The reduced normal equation for $\widehat{\alpha^T}$  is given by
 \begin{equation}  \label{ReducedNEalphaT}
 {\mathbf C}_{T;\bar{T}} \widehat{\alpha^T} = {\mathbf Q}_{T;\bar{T}}, \mbox{ where }
 \bar{T} = {\cal F} \setminus T.
 \end{equation}

 In particular, the reduced normal equation for  $\widehat{\alpha^A}$  is
\begin{equation}  \label{ReducedNEalphai}
 {\mathbf C}_A \widehat{\alpha^A} = {\mathbf Q}_A.
 \end{equation}

(c) The sums of squares considered in (e) of Notation \ref
{C-matandQ} can be expressed as follows.
\begin{eqnarray*} SS_{A;T} & = & {\mathbf Q}'_{A;T} ( C_{AA;T})^- Q_{A;T} \\
\mbox{ and }  SS_{S;T} & = &{\mathbf Q}'_{S;T} ({\mathbf C}_{S;T})^-
{\mathbf Q}_{S;T}
\end{eqnarray*}

In particular  the sum of squares for $A$, adjusted for all the
other factors is $SS_{A;\bar{A}} = {\mathbf Q}'_A ({\mathbf C}_A)^{-} {\mathbf Q}_A$.
\end{lem}

\vspace{.5em}

We need the following well-known results, which also follows from
Lemma \ref {NEandSS}

\begin{lem} \label{ssadj}
(a) For $A \neq B,\; A,B \in {\cal F}$,  $SS_{A;B}$ is the quadratic form $Y'P_U Y$, where $U
=(I-P_{B})X_A$.

(b) More generally, for $S,T \subset {\cal F}, \;S \cap T = \phi$, $SS_{S;T}$ is the quadratic form $Y'P_V Y$, where $V =(I-P_T) X_S$.

In particular, $SS_{A;\bar{A}} = Y'P_H Y$, where $H =(I-P_{\bar{A}}) X_A$.

(c) The so-called unadjusted sum of squares for $A$ is $SS_{A;G} = A' (R_A)^{-1} A -G^2/n$, where $A$ is the vector of raw totals
for $A$ and $G$ is the grand total.

\end{lem}

\vspace{.5em}

{\bf Plans for symmetric experiments laid out in blocks.}

We now assume that a block factor is present.  However, it is convenient to view the block factor separately from the treatment factors. We, therefore, use the following notation.

\begin{nota} \label{exptBl}
 (a)   $b$  will denote the number of blocks and $k_j$  the
 size of the $j$th block, $\;1 \leq j \leq b$. Thus, $b = S_B$. Also,  the total number
 of runs $n = \sum_{j=1}^{b} k_j$.  $\{ B_j, \: j=1, \cdots b\}$
 will denote the set of all blocks of ${\cal P}$.

(b) $ {\cal T}$ will denote the set ${\cal F} \setminus \{B\}$, where $B$ is the block factor. $L_A$ will denote the $A$-versus block incidence matrix, $A \in {\cal T}$. Thus, $L_A = N_{AB}$ and the $(p,j)$th entry of  $L_A$ is
$$l^A (p,j) = |\{x \in B_j :  x_A = p\} |, \; p \in S_A, 1 \leq j \leq b,\;A \in {\cal T}.$$

(c) $D_k$ will denote the diagonal matrix whose diagonal entries are  $k_j, \;1 \leq j \leq b$, in that order. Thus, $D_k = R_B$.
\end{nota}

{\bf Reduced normal equation for the contrasts :}
Sometimes it is useful to consider the reduced normal equation for the contrasts
and the corresponding C-matrix. Towards that we introduce the following notations.

\begin{nota} \label{Contrast}
(a) For each   factor $A$,  $O_A$ will denote an $(s_A - 1) \times s_A$ matrix such that $O_A O'_A = I_{s_A-1}$ and $O_A 1_{s_A} = 0$.
Let $Z_A = X_A O'_A, \; \gamma^A = O_A \alpha^A$.
\end{nota}
Then, the model (\ref {modelEq}) can also be  expressed  as
\begin{equation}\label{modelContrast} {\mathbf Y} =  {\mathbf X}_G  \alpha^G + {\mathbf Z} \delta +\epsilon
\mbox{ where } {\mathbf Z} = [ {\mathbf Z}_A \; : \; A \in {\cal F} \setminus\{G\}] \mbox{ and } \delta = [\gamma^A \; : \; A \in {\cal F} \setminus\{G\}]'. \end{equation}
Here  ${\mathbf Z}_A$'s and $\gamma^A$'s are as in Notation \ref {Contrast},
while ${\mathbf X}_G,  \alpha^G$ and $\epsilon$ are as in Notation \ref {modelEq}.

It is easy to verify that
\begin{equation}\label{NEcontrast}
{\mathbf Z}' {\mathbf Z} \widehat{\delta} = {\mathbf Z}' {\mathbf Y}. \end{equation}

Now we assume that a block factor is also present. We obtain the reduced normal equation
for the vector of all contrasts, after eliminating the block effects.
Towards this, we use Notation \ref {exptBl} and  Notation \ref {C-matandQ} (c) with $T = {\cal T}$
and $U = \{B\}$.

\begin{eqnarray}\label{NEcontrastBl}
 \tilde{C}  \widehat{\delta} =  Q^P, \mbox{ where } & \tilde{C}  =  (( \tilde{C}_{AA'}))_{A,A' \in {\cal T}}, & \tilde{C}_{AA'} =
 O_A (N_{AA'} - L_A D^{-1}_k L'_{A'}) O'_{A'} \\
 \mbox{and } & Q^P =  ((Q^P_A))_{\cal T}, &  Q^P_A = O_A X'_A (I - P_B) Y.
 \end{eqnarray}

\begin{defi}\label{C-contrast} By the C-matrix of the contrasts of a plan ${\cal P}$, we shall mean a matrix $\tilde{C} ({\cal P})$. This is  the matrix ${\mathbf Z}' {\mathbf Z}$ of (\ref {NEcontrast}) if no block factor is present, while it is $\tilde{C}$ in (\ref {NEcontrastBl})
if a block factor is present.
\end{defi}

{\bf Remark :} The   C-matrix of the contrasts for a plan ${\cal P}$ is particularly useful when partially orthogonality holds among one or more factors of ${\cal P}$ [see Definition 2.4 of Bagchi (2019)]. As an example let us take ${\cal P}$ to be the plan $ICA (N,3^l2^p)$ of Huang, Wu, and Yen, C.H. (2002), where $p = N - 2(l+1)$. One can see that $\tilde{C} ({\cal P})$ is of the form $\tilde{C} ({\cal P})  =
\left [   \begin{array}{ccc} C_1 & 0 \\
0  & C_2 \end{array} \right ]$,
where $C_2 = d I_p$ for a real number $d$. Further, conjugation $C_1$ by a suitable permutation matrix,  we can write it as $\left [   \begin{array}{ccc} A_1 & 0 \\
0  & A_2 \end{array} \right ]$, $A_i = a_i I_l + b_i J_l, i = 1,2$. Thus,  one can see that the set of all contrasts of all the three-level factors satisfy inter-class orthogonality, the orthogonal classes being the  linear contrasts and the quadratic contrasts.

\section{Orthogonality through other factors  versus usual orthogonality}
\setcounter{equation}{0}
In this section we seek the answer to the following questions. {\bf Consider a
main effect plan for $m (\geq 3)$ factors. Fix a factor, say $A$.
What conditions  must the design matrices satisfy so that
the inference on $A$ depends only on the relation of $A$ with the factors in a certain class of factors
(say $T$) ? }
In other words, for inference on $A$ one may forget all factors other than those in $T$.

 Towards an attempt  to answer these questions, we need a definition.

\begin{defi}\label{orth-S} (a) Consider $T \subset  ${\cal F}$, \; A,B \notin \ T$.
  Then, the factors $A$ and $B$ are said to be
 {\bf  orthogonal through the factors in $T$ }, denoted by $A \bot_T B$, if
\begin{equation}  \label{orthogonal-S}
  N_{AB} = N_{AT} (X'_T X_T)^{-1} N'_{BT} \mbox{ or equivalently, } X_A (I - P_T) X_{B} = 0.
\end{equation}

(Here $N_{AT}$ is as in   (i) of Notation \ref {C-matandQ} (b)).

\end{defi}

We need the following well-known results. [See Exercise 7 of Chapter 2 of  Yanai, H., Takeuchi, K. and  Takane, K. (2011), for instance].

\begin{lem} \label{PA-PD} Consider matrices $U,V,W$ with the same number of rows. Then, the following hold.

(a) Suppose $ {\cal C} (V) \subseteq  {\cal C} (W)$.  Then
$$ {\cal C} (P_V U) = {\cal C} (P_W U) \Leftrightarrow (P_W - P_V)U =0.$$
(b) If $W = [U,V]$ then, $ P_W - P_V = P_Z$, where $Z = (I - P_V) U$.
 \end{lem}



We are now in a position to seek answers to the question posed above. We shall use the abbreviation w.p.1 for the phrase with probability
1.

\begin{theo} \label{EstSSQ}
Fix $A \in {\cal F}$. Partition $\bar{A}$ as $\bar{A} = S \cup T$. Then, a
necessary and sufficient condition for each of the following statements is that
$A \bot_T B, \forall B \in S$.

(a) $C_{A;T}$ is the C-matrix for $A$.

(b) $SS_{A;\bar{A}} = SS_{A;T}$ w.p.1
\end{theo}

{\bf Proof :} Taking $W = X_{\bar{A}}, \; U = X_S$ and $V = X_T$  and
applying Lemma \ref  {PA-PD} we find that
\begin{equation}  \label{Pmatrices} P_{\bar{A}} =  P_T  + P_Z, \mbox{ where } Z = (I - P_T) X_S. \end{equation}
{\bf Proof of (a):} In view of (\ref {CQSS}), (\ref {ReducedNEalphai})and  (\ref {Pmatrices})
we see that $C_A = C_{A;T} - X'_A P_Z X_A$. Therefore, the required necessary and sufficient condition
is that $P_Z X_A = 0$, which is equivalent to (\ref {orthogonal-S}). Hence the result.

{\bf Proof of (b):}  In view of Lemma \ref {ssadj} the following hold.
\begin{eqnarray}  \label{sumSq}
SS_{A;\bar{A}} = Y^{\prime }P_U Y, & SS_{A;T} = Y^{\prime }P_V Y,\\
\mbox{ where } U = (I - P_{\bar{A}}) X_A, & V = (I - P_T)X_A.\end{eqnarray}

Since the support of $Y$ is $R^n$ (recall (\ref {modelEq})),
$Y^{\prime }P_U Y = Y^{\prime }P_V Y$  w.p.1 if and only if $P_U =
P_V$. Thus, the required necessary and sufficient condition is that
${\cal C} (U) = {\cal C} (V)$. Applying  Lemma \ref {PA-PD} we see that the required necessary and sufficient
condition is
\begin{equation}  \label{cond}
( P_{\bar{A}} - P_T) X_A = 0.
\end{equation}
But in view of (\ref {Pmatrices}) this is the same as  $P_Z X_A = 0$, which is equivalent to (\ref {orthogonal-S}). Hence the result. $\Box$

Extending Definition \ref {orth-S}, we have the following.

\begin{defi} \label{planOrthT} Consider a plan ${\cal P}$ with $m$ factors.  Suppose there is a $T \subset {\cal F}$  of size $t$  such that
$A \bot_T B$ for every pair $(A,B)$ such that $A \neq B$ and $A,B \notin T$. Then we say that  ${\cal P}$ is a {\bf plan
orthogonal through $t$ factors}. \end{defi}

{\bf Special cases :}

{\bf Case  $|T| = 1$ and $T = \{G\}$ :} This reduces Definition \ref {orth-S} to the usual definition of orthogonality and (\ref {orthogonal-S}) to the  proportional frequency condition (PFC) of Addelman (1962). We present this important condition below. Two factors $A \neq B$ are orthogonal if
\begin{equation} \label{PFC} n N_{AB} = r_A r'_{B} .
\end{equation}


This case of  Theorem  \ref {EstSSQ}  is   well-known.

{\bf Case  $|T| = 1$, but  $T \neq G$ :} This case is referred to as  orthogonality between a pair of factors through a third factor in Bagchi (2019), where $T$ is a treatment factor. Examples of plans satisfying the conditions of Definition \ref {planOrthT} may be found in Section 2 of the same paper. The case when $T = \{B\}$, i.e. orthogonality through a block factor is considered in Bagchi (2010). In that paper a plan of Definition \ref {planOrthT} with $T = \{B\}$ is termed a plan orthogonal through the block factor (POTB). By now many POTBs are available in the literature. Reference :

{\bf Case $|T |= 2$ :}  In this case the plan may be termed as a  plan orthogonal through a pair of factors. We shall construct a series of such plans in Section 4.1.

{\bf Optimality Criteria :}
Among the plans available in a given set up, one would like to use the one which estimates the contrasts of interest more precisely than the other plans. To assess the performance of a plan, one applies  an ``optimality criteria"
on its `information matrix" or ``C-matrix". Now, regarding comparison between plans, one may look at the plan as a whole or look at its performance regarding one or more factors.  We may note here that the first approach is particularly meaningful when the plans are compared in terms of the C-matrix of all the contrasts. In the second approach, one chooses a factor $A$ and then compare the plans  in terms of $C_A$ (see (d) of Notation \ref {C-matandQ}), using an optimality criterion. This  approach has been used in Bose and  Bagchi (2007), where a plan was found to be E-optimal for two factors, while universally optimal for the other two. For the definition of universal optimality and other details we refer to Shah and Sinha (1989).

Modifying the celebrated  Theorem of Kiefer (1975)  to be applicable to the  context of plans with a block factor,  we get the following.

\begin{theo}\label{univOpt} Consider a class $\Pi$ of connected plans   as described in  Notation \ref {exptBl} and a plan ${\cal P}^* \in \Pi$.
 Let $t_j = [k_j/s_A], 1 \leq j \leq b, A \in {\cal T}$. Consider  the following conditions.

(a) For a fixed $A, A \in {\cal T}$,  ${\cal P}^*$ satisfies the following.

(i) In each of the $s_A$ levels of $A$ appears $t_j$ or $t_j + 1$ times in the $j$th block, $1 \leq j \leq b$.

(ii) $A \bot_B A', \; A' \neq A  \in {\cal T}$ in ${\cal P}^*$ [recall that $B$ is the block factor]

(iii) $C_A ({\cal P}^*)$ is of the form $aI + b J$.

\vskip5pt

(b) $\tilde{C} ({\cal P}^*)$ is of the form $aI_v$, where $v = \sum \limits_{A \in {\cal T}} (s_A-1)$.

We have the following.

(a) If a plan ${\cal P}^* \in \Pi$ satisfies all the conditions in (a), then ${\cal P}^* \in \Pi$ is universally optimal in $\Pi$
for the inference on $A$.

(b) The pair of conditions (a)(i) and (b) is sufficient for  universally optimality of ${\cal P}^*$  for the inference on all the main effect contrasts.
\end{theo}

We present the definition of an widely used optimality criterion, which we need in the next section. For a real symmetric $n \times n$ matrix $A$, $ \mu_0 (A) \leq \cdots \leq \mu_{n-1} (A)$ will denote the
eigenvalues of $A$.

\begin{defi} \label{ADEopt}

  A plan ${\cal P}^* \in \Pi$ is said to be E-optimal in for the inference on  all the main effect contrasts
      if $ \mu_1 (\tilde{C} ({\cal P}^*)  \geq \mu_1 ( \tilde{C} ({\cal P}) ) \; \; \forall {\cal P} \in \Pi$.
 \end{defi}

\section{Construction of plans orthogonal through one or two factors}
\setcounter{equation}{0}

We shall now proceed to construct main effect plans,  mostly  for  symmetric experiments. In Section 4.1
the factors are assumed to be treatment factors, while in Section 4.2 a block factor is assumed to be present.
 Most of the constructions are of recursive type, in the sense that from a given initial plan we generate a plan by adding blocks and/or factors.

\begin{nota}\label{AddBl}
(a) $s$ is a prime power. $F$ will denote the Galois field of order $s$. $F^m$ will denote the vector space of dimension $m$ over $F$.

(b) The set of levels of each factor of a plan  is $F$ unless stated otherwise.
 For every plan ${\cal P}$,  $R({\cal P})$ will denote the set of all runs of ${\cal P}$. Thus, if ${\cal P}$ is a plan for an $s^m$ experiment, then $R({\cal P}) \subset F^m$.

(c) Consider an initial plan ${\cal P}_0$  for an $s^m$ experiment. Consider $V \subset F^m$ (a set of runs for the same experiment as ${\cal P}_0$).  ${\cal P}_0 + V$ will denote the plan having
the set of all runs $\{x + v, v \in V, x \in  R({\cal P})\}$. Here the addition is the vector addition in $F^m$.

In particular, if $V = \{u1_m, u \in F\}$, then ${\cal P}_0 + V$ will
be denoted simply by ${\cal P} \oplus F$.
\end{nota}

\begin{defi} \label{AddingBlock}
Consider an initial plan ${\cal P}_0$ for an $s^m$ experiment on $n$ runs, which are  distinct members of $F^m$.
Consider a subset $V$ of $S^m$ containing the $0$ vector. By the plan {\bf generated from ${\cal P}_0$ along $V$} we shall mean the plan (for the same experiment) ${\cal P}_0 + V$ as in (c) of Notation \ref {AddBl}. The set $V$ is termed {\bf generator}.
\end{defi}

\subsection{ Plans orthogonal through a pair of factors.}

We shall construct a series of plans orthogonal through a pair of factors (POTP). The condition for such an orthogonality is rather strong and
so a POTP is rather rare. The following result provides a direction for the search of a POTP.

\begin{theo} \label{orthogonal-2} Consider an $s^m$ experiment.
 If the incidence matrices are as shown below, then
$A_i$ and $A_j$ are orthogonal through the pair $(A_1,A_2)$, for every pair $(i,j), \; i,j > 2$. Here  $c$ is an integer.

\begin{equation}  \label{condition2-orth}
  N_{12} = 2c(J_s - I_s), \mbox{ and } N_{ij} = c((s-2)I_s + J_s), \; i \neq j, \; (i,j) \neq (1,2).
\end{equation}
\end{theo}

The proof is  by  straightforward verification of (\ref {orthogonal-S})
with $T = \{1,2\}$.

\vspace{.5em}

{\bf Example :} The following is a plan for a $3^4$ experiment  satisfying the conditions of Theorem \ref
{condition2-orth}.

\begin{center}
{\bf {\large Table 4.1.1} :}

\vspace{.5em}

\begin{tabular}{ll|cccc cccc cccc}
 \hline
  Factors $\downarrow$ &  $A_1$  & 0 & 0 & 0 & 0 & 1 & 1 & 1 & 1 & 2 & 2 & 2 & 2 \\
                       &  $A_2$  & 1 & 1 & 2 & 2 & 2 & 2 & 0 & 0 & 0 & 0 & 1 & 1 \\
                        & $A_3$  & 0 & 1 & 2 & 0 & 1 & 2 & 0 & 1 & 2 & 0 & 1 & 2 \\
                        & $A_4$  & 0 & 1 & 0  &2 & 1 & 2  & 1 & 0 & 2 & 0 & 2 &  1\\
                     \hline
            \end{tabular}
\vskip5pt
\end{center}
We shall now present a general  construction for a POTP for a symmetric experiment.

\begin{theo} \label{2orthPlan} Suppose $h \equiv 0 \pmod 4$ and $h$ is the order of a Hadamard matrix. If $s \equiv 3 \pmod 4$ is a prime power, then there exists a POTP for a $s^{2h}$ experiment on $2hs(s-1)/2$ runs.
\end{theo}

Towards the proof of this theorem we define the following.

\begin{nota}\label{sPrime}
 (a) $s \equiv 3 \pmod 4$ is a prime power. $F^* = F \setminus \{0\}$. $C_0$   will denote the set of all non-zero squares of $F$, while $C_1$ will denote the set of all non-zero non-squares of $F$.

 (b) Consider  an $m \times n$ array $P$ with entries from $F$.  For $1 \leq i,j\leq m, \; k \in F$, let
  $$d^k_{ij} = |\{l : p_{jl} - p_{il} = k \}|.$$
   $C_0 P$ will denote the $m \times (s-1)n/2$ array  $ \{cP: c \in C_0\}$.
  \end{nota}

\begin{lem}\label{Nij} Suppose $\exists$ an  array $P$ as in Notation \ref {sPrime} satisfying the following.
$\sum \limits_{ k \in C_0} d^k_{ij} = \sum \limits_{ k \in C_1} d^k_{ij} = u_{ij}$,  (say). Let ${\cal P}$ be a  plan having the columns of the $m \times (s-1)n/2$ array $C_0 P$ as the set of runs. Let ${\cal P}^* = {\cal P} \oplus F$. Then, for a pair $(i,j), i\neq j,\;1 \leq i,j \leq m$,  the incidence matrix $N_{ij}$ of ${\cal P}^*$ satisfies the following.
$$N_{ij} (x,y) =  w _{ij}I_s + u_{ij} J_s, \mbox{ where } w _{ij} = (s-1)n/2 - su_{ij}. $$
\end{lem}

{\bf Proof :} $N_{ij} (x,y) = |\{(l,\alpha,q) : 1 \leq l \leq n, \alpha \in F, \;q \in C_0, y-x = q(p_{jl} - p_{il}), \alpha = x - qp_{il}\}|$, which is $=|\{(l,q) : 1 \leq l \leq n, q \in C_0, y-x = q(p_{jl} - p_{il})\}|$. So, by hypothesis,
$$\mbox{ if $y-x \in C_0$, then }
 N_{ij} (x,y) = |\{ l : p_{jl} - p_{il} \in C_0\}| = \sum \limits_{ k \in C_0} d^k_{ij} = u_{ij}.$$
 Similarly, if $y-x \in C_1$, then $ N_{ij} (x,y) = \sum \limits_{ k \in C_1} d^k_{ij}$.  Therefore, if $y = x$, then
$ N_{ij} (x,y) =  \frac{s-1}{2}(n - 2 u_{ij})$. Hence the result. $\Box$

Using the fact that when $ s  \equiv 3 \pmod 4$, $-1 \in C_1$, we get the following result from Lemma \ref {Nij}.

\begin{lem}\label{1-1 array} Suppose $ s  \equiv 3 \pmod 4$  is a prime power and  $n$ is a multiple of $4$. Suppose there is an
$m \times n$ array $P$ with entries $0,1,-1$ (viewed as members of $F$) satisfying the following.

(a) $p_{1,j} = 0, 1 \leq j \leq n$.

(b) For every ordered pair $(i,j)$,  $p_{il} - p_{jl} \in \{0,1,-1\}, \; 1 \leq l \leq n$.

(c) For $k = 1,-1, d^k_{ij} = \left\{ \begin{array}{ll} n/2 & \mbox{ if } (i,j) = (1,2), \\
n/4 & \mbox{otherwise} \end{array} \right. $

Then the plan ${\cal P}^* = C_0 {\cal P} \oplus F$ satisfies the conditions of Theorem
\ref {orthogonal-2} with $c = (n/4)s(s-1)/2 = ns(s-1)/8$.
\end{lem}

{\bf Proof of Theorem \ref {2orthPlan} : } Let $H$ be a Hadamard matrix of order $h$. W.l.g., we assume that the
first row of $H$ consists of only $1$'s.
Write $H = \left [ \begin{array}{cc} 1_h & \tilde{H} \end{array} \right ]' $. Consider the array
$$P = \left [\begin{array}{ccc} 0_{1 \times h} &|&  0_{1 \times h} \\
 J_{1 \times h} &|&  - J_{1 \times h} \\
 (\tilde{H} + J_{h-1 \times h} )/2 &|&  - (\tilde{H} + J_{h-1 \times h} )/2\\
(\tilde{H} + J_{h-1 \times h} )/2 &|&   (\tilde{H} - J_{h-1 \times h} )/2\\
 \end{array} \right]. $$

It is easy to check that $P$ satisfies the conditions of Lemma \ref  {1-1 array} and hence the result follows
from the Theorem \ref {orthogonal-2}. $\Box$

\subsection{ Plans orthogonal through the block factor}
In this section we construct main effect plans plans orthogonal through the block factor. We, therefore need the definition of orthogonality through the block factor. This is obtained by taking $T$ to be a block factor in (b) of  Definition \ref {orth-S}. However, it will be convenient to express the Definition \ref {orth-S} using quantities in Notation \ref {exptBl}. It is interesting to note that  Condition (\ref {orthogonal-bl}) below  is equivalent to equation (7) of Morgan and Uddin (1996) in the context of nested
row-column designs.

\begin{defi} \label{Orth-bl}
(a)  Fix $A \neq A', \: A,A' \in {\cal T}$. The factors $A$ and $A'$ are said to be
 {\bf  orthogonal through the block factor } if
\begin{equation}  \label{orthogonal-bl}
  N_{AA'} = L_A (D_k)^{-1} L_{A'}.
\end{equation}


(b) A plan ${\cal P}$ is said to be a {\bf plan orthogonal through the
block factor (POTB)}  if (treatment) every factor is orthogonal to every other one through the block factor.

\end{defi}

{\bf Remark :} Now onwards  orthogonality will mean orthogonality through the block factor.

  In a recursive construction we need to find a suitable initial plan ${\cal P}_0$ as well as a suitable generator $V$, so that the generated plan ${\cal P}_0 + V$ satisfies one or more desirable property. We now look for such a suitable generator.

{\begin{nota} Consider $V \subset F^m$, viewed as a set of runs for an $s^m$ experiment with ${\cal T}$ as the set of factors. Fix $v \in V$.
 For $ A \in {\cal T}, \; v_A$ is the level of $A$ in the run $v$. Similarly, for $A,A' \in {\cal T}, v_{AA'} $   will denote the ordered pair
$(a,a')$, where $a$ (respectively $a'$) is the level of $A$ (respectively $A'$) in $v$.  $V_A$ and $V_{AA'}$ will denote the following multisets of size $|V|$  contained in $F$ and $F \times F$ respectively.
$$V_A = \{v_A : v \in V\}, \mbox{ while } V_{AA'} = \{ v_{AA'} : v \in V \}.$$
\end{nota}

The following  useful result can be verified easily.

\begin{lem}\label{AiAjOrth} If every element of $F \times F$ appears a constant number of times in $V_{AA'}$,
then, $A \bot A' $ in ${\cal P}_0 + V$, no matter what ${\cal P}_0$ is.
\end{lem}

We shall now see how we can  enlarge the set of factors of a given plan, while keeping the
number of blocks fixed.

\begin{defi}\label{powerPlan} (a) Consider a plan ${\cal P}$ as in Notation \ref {exptBl} and an integer $t$. Let $x_{ij}$ denote the $j$th run in the ith block of ${\cal P}$. Consider the $t(m-1) \times 1$ array $y_{ij}$ obtained by juxtaposing the $(m-1) \times 1$ array $x_{ij}$ $t$ times,
Then, the plan on b blocks  with $y_{ij}$ as the $j$th run in the $i$th block, $1 \leq j \leq k_i, \: 1 \leq i \leq b$ is said to be obtained by taking the $t$th power of  ${\cal P}$. The new plan will be denoted by ${\cal P}^t$.  We name the factors of ${\cal P}$ and its power ${\cal P}^t$ as in Notation \ref {allFactors} below.
\end{defi}

\begin{nota}\label{allFactors} Consider a plan  ${\cal P}$ having
a set of factors ${\cal T}_0 = \{A, \cdots M\}$. The set of
factors of ${\cal P}^t$ will be named as
$${\cal T} = \bigcup\limits_{i = 1}^t {\cal T}_i, \mbox{where }
 {\cal T}_i = \{A_i, \cdots M_i \}.$$

\end{nota}


Combining Definitions \ref {AddingBlock} and \ref {powerPlan} we get a recursive construction in which factors as well as blocks are added to the initial plan.

\begin{defi} \label{recursiveConstr}
 Consider an initial plan ${\cal P}_0$  for an $s^m$ experiment laid in $b$ blocks
 of sizes $k_1, \cdots k_b$. Consider a $p \times q$ array $H = ((h_{ij}))_{1 \leq i \leq p, 1 \leq j \leq q}, \; h_{ij} \in F$.  We now obtain a plan for an $s^{mq}$ experiment on $bp$ blocks using the array $H$ as follows. We first  obtain ${\cal P}_0^q$ following Definition
\ref {powerPlan}.

Let $w_i = \left [ \begin{array}{ccccccc} h_{i1}.1^{\prime }_t &
h_{i2}.1^{\prime }_t & \cdots & h_{iq}.1^{\prime }_t \end{array}
\right ]^{\prime }, \: 1 \leq i \leq p$ and $W_H = \{w_i, \: 1 \leq
i \leq p\}$.

Our required plan ${\cal P}$  is $({\cal P}^q_0) + W_H$ and it will be denoted by $H \Diamond {\cal P}_0$. Symbolically,
\begin{equation}  \label{NewPlan}
{\cal P} = H \Diamond {\cal P}_0 = {\cal P}_0^q + W_H.
\end{equation} \end{defi}

Here the  factors of ${\cal P}_0^q$ as well as ${\cal P}$ are named in accordance with  Notation \ref {allFactors}.

 Our task is to find a suitable array $H$ so that the plan $H \Diamond {\cal P}_0$ satisfies certain desirable properties. The
natural choice would be  an orthogonal array of strength two [see Rao(1946)]. We shall use a
  slightly modified version of it, so as to accommodate a few more factors.

\begin{nota}\label{ExtOA}  An orthogonal array of $m$ rows, $N$ columns with the entries  from the set of integers modulo $s$ and  strength $t$ will be denote by $OA(N,m, s, t)$.

 The array obtained by adding a row of all zeros (in the $0$th
position, say) to an $OA(N,m-1,s,2)$ will be denoted by $Q(N,m,s)$.
\end{nota}



\vspace{.5em}

 We get the following result from the recursive construction described
in Definition \ref {recursiveConstr}.

\begin{theo} \label{PlanThruOA} Consider a plan ${\cal P}_0$ for an $s^t$ experiment on $b$ blocks. Suppose an orthogonal array
$OA(N,m-1,s,2)$exists. Then $\exists$ a plan ${\cal P}$ with a set of $s^{mt}$ factors on $bN$ blocks with the following properties. Here the names of the factors are in accordance with Notation \ref {allFactors}.

(a) If ${\cal P}_0$ has $b_j$ blocks of size $k_j$, then ${\cal P}_N$
has $Nb_j$ blocks of size $k_j, \; j = 1, \cdots b$.

(b) Consider $P \neq Q, P,Q \in {\cal T}_0$. If  $P \bot Q $ in ${\cal
P}_0$, then $P_i \bot Q_i $ for every $i,\: 0 \leq i \leq m - 1$,

(c) $P_i \bot Q_j , \; P,Q \in {\cal T}_0,\; i \neq j, 0 \leq i,j \leq
m , $. \end{theo}

{\bf Proof :} By hypothesis, the array $Q(N,m,s)$ exists. The required plan  ${\cal P}$ is $Q(N,m,s) \Diamond {\cal P}_0$. Properties  (a) and (b) follow from the construction while (c) follows from Lemma \ref {AiAjOrth}.

\vspace{.5em}

{\bf New POTB for two-level factors.}

We now present an initial plan ${\cal P}_0$ for seven factors on two
blocks $B_1$ and $B_2$, each of size five.
 \begin{center}
{\bf  Table 4.2.1 : The initial plan $ {\cal P}_0$}

\vskip5pt

\begin{tabular}{ll|ccccc|ccccc}
Blocks &$ \rightarrow$ & &  & $B_1$ &  & & &  &$B_2$ & &  \\
 \hline
  Factors $\downarrow$ &  $A_1$  & 0  &0 & 1 &  1 & 0 & 0 & 1 & 1 & 1 &  1  \\
                       &  $A_2$  & 0  & 1 & 0 & 1 & 0 & 1 & 0 & 1 & 1 &  1 \\
                        & $A_3$  & 0  & 1 & 1 & 0 & 0 & 1 & 1 & 0 & 1 &  1 \\
                        & $A_4$  & 0  &0 & 1 & 1  & 1 & 1 & 0 & 0 & 0 &  0\\
                        & $A_5$  & 0  &1  & 0& 1  & 1 & 0 & 1& 0 &0 &0 \\
                        & $A_6$  &0  &1 & 1 & 0 & 1 & 0 & 0 & 1 & 0 &0 \\
                         & $A_7$ &0  &0 & 0 & 0 & 1 & 0 & 0 & 0 & 1 &1  \\\hline
            \end{tabular}
 \end{center}
\vspace{.5em}

\begin{theo}\label{2^7} (a) The plan ${\cal P}_0$ is a  POTB for a $2^7$
experiment.

(b) Moreover, it is E-optimal among all plans for the same
experiment on the same set up.
\end{theo}

{\bf Proof : } (a)  is proved by straightforward verification of (\ref {orthogonal-bl}) for each pair $(i,j) : i \neq j,\; i,j = 1, \cdots 7$.

(b) Using (\ref {NEcontrastBl}), the C-matrix of the contrasts of ${\cal P}_0$  is obtained as
\begin{equation}\label{MP0}
\tilde{C} ({\cal P}_0) = 4I_7.\end{equation} The rest follows from Theorem 3.1 of Jacroux and Kealy-Dichone (2015). [ E-optimality is defined in definition \ref {ADEopt}] $\Box$

\vspace{.5em}

{\bf Remark 5.2:} The set up of ${\cal P}_0$ is an example of Case 3 of Jacroux and Kealy-Dichone (2015). The E-optimal plan
constructed in that paper in the same set up has at most four factors, while ${\cal P}_0$ accommodates seven factors. The model used in Jacroux and Kealy-Dichone (2015) is a bit different from  the one given in Notation \ref {Contrast}. This is due to the fact that the rows of the matrix $O_A$ in their model are not orthonormal. If we use their model, then $ \tilde{C} ({\cal P}_0)$ would be $8I_7$.

\vspace{.5em}

Next we derive an infinite series of POTBs from ${\cal P}_0$.

\begin{theo}\label{2power7n} (a) If there exists a  Hadamard matrix of order $h$ $\geq 2$, then
 there is a  POTB ${\cal P}_h$ for a $2^{7h}$ experiment in $2h$ blocks of size $5$ each.

(b) ${\cal P}_h$  is E-optimal among all plans for the same
experiment on the same set up.

\end{theo}

{\bf Proof :} Suppose $h \geq 2$ is a  Hadamard number. By the well-known relation between orthogonal arrays of strength two and Hadamard
matrices, (see Theorem 7.5 in Hedayat, Sloane and Stuffken (1999), for instance), an $OA(h-1,h,2,2)$ and hence  $Q(h,h,2)$ exists.

 (a)  Let ${\cal P}_h = Q(h,h,2) \Diamond {\cal P}_0$, where ${\cal P}_0$ is as in Table 5.1 (recall
 Definition \ref {recursiveConstr}). Now, Theorem \ref {2^7} and Theorem \ref {PlanThruOA} together
  imply that the plan ${\cal P}_h$ is a POTB.

(b) Theorem \ref {PlanThruOA} together with \ref {2^7}
 imply that the information matrix of ${\cal P}_h$ is $4hI_{7h}$.
   Hence, by the same argument as in Theorem \ref {2^7} the E-optimality of
${\cal P}_h$ follows. $\Box$

\vspace{.5em}

We construct anther plan in the same set up as in Theorem \ref {2^7}, but with $6$ factors.

\begin{center}
{\bf  Table 4.2.2 :  Plan $ {\cal P}_2$}
\vskip5pt

\begin{tabular}{ll|ccccc|ccccc}
Blocks &$ \rightarrow$ & &  & $B_1$ &  & & &  &$B_2$ & &  \\
 \hline
  Factors $\downarrow$ &  $A_1$  & 0  & 0 & 1 & 1 & 0 & 0 & 0 & 1 & 1 &  0  \\
                       &  $B_1$  & 0  & 1 & 0 & 1 & 0 & 0  & 1 & 0 & 1 & 0 \\
                        & $C_1$  & 0  & 1 & 1 & 0 & 0 &  0  & 1 & 1 & 0 & 0 \\
                        & $A_2$  & 0  & 0 & 1 & 1 & 0 & 1 & 1 & 0 & 0 &1 \\
                        & $B_2$  & 0  & 1 &  0 & 1& 0 & 1 & 0 & 1 & 0 & 1\\
                        & $C_2$  & 0  & 1 & 1 & 0 & 0 & 1& 0 & 0 & 1 &1  \\
            \end{tabular}
 \end{center}
\vspace{.5em}

This is an inter-class orthogonal plan, the classes being $\{A_1, B_1, C_1, D_1\}$ and $\{A_2, B_2, C_2, D_2\}$.
The  C-matrix of the contrasts $\tilde{C} ({\cal P}_2) = \left [ \begin{array}{ccc} A & 0 \\
                                                                              0 & A  \end{array} \right ]$, where $A = 8(I_3 + (1/5) J_3$.

{\bf Discussion :}  Consider  the set  up of in Theorem \ref {2^7} with six factors. The hypothetical universally optimal plan, say ${\cal P}^*$ would have  C-matrix of the contrasts $\tilde{C} ({\cal P}^*) = (48/5)I_6$. Let $\tilde{{\cal P}}_0$ denote the plan obtained by deleting the factor $A_7$ from  $\tilde{{\cal P}}_0$. We see that both ${\cal P}_2$ and $\tilde{{\cal P}}_0$ are E-optimal, but ${\cal P}_2$ is A-better that  $\tilde{{\cal P}}_0$.

Again, we note that the plan $\tilde{{\cal P}}_0$ satisfies Conditions  (b) of Theorem \ref {univOpt}, but not (a) (i). On the other hand,
${\cal P}_2$ satisfies Condition (a)(i), but do not satisfy (b). [Recall that by  Theorem \ref {univOpt}, the pair of conditions (a)(i) and (b)  is sufficient for universal optimality.]
This is the usual story in the search for optimal designs, which is why an universally optimal design is rare.  We conjecture that the hypothetical  plan ${\cal P}^*$ does not exist and ${\cal P}_2$ is A-optimal.

 \vspace{.5em}

{\bf  New POTB for three-level factors}

The following is an initial plan ${\cal P}_0$ for three factors on three
blocks   of sizes $4,4,2$.
 \begin{center}
{\bf  Table 3.2.3 : The initial plan $ {\cal P}_0$}

\vskip5pt

\begin{tabular}{ll|cccc|cccc|cc}
Blocks &$ \rightarrow$ & &  & $B_1$ &  & &   &$B_2$ && $B_3$  \\
 \hline
  Factors $\downarrow$ &  $A_1$  & 0  &0 & 1 &  2 & 0 & 0 & 1 & 2 & 1 &  2  \\
                       &  $A_2$  & 0  & 1 & 0 & 2 & 2 & 0 & 1 & 0 & 2 &  1 \\
                        & $A_3$  & 0  & 1 & 2 & 0 & 2 & 0 & 0 & 1 & 1 & 2 \\
                      \hline
            \end{tabular}
 \end{center}
\vspace{.5em}

\begin{theo}\label{3^3} (a) The plan ${\cal P}_0$ is a  POTB for a $3^3$
experiment.

(b) Moreover, ${\cal P}_0$ is universally optimal among all plans for the same
experiment on the same set up.
\end{theo}

{\bf Proof :} (a) follows from the fact that  (\ref {orthogonal-bl}) holds for each pair
$(A,A') : A,A' \in {\cal T}$.

To prove (b) we note that the condition (a) (i) of Theorem \ref {univOpt} is satisfied. Then, we
 use (\ref {NEcontrastBl}) and get that
\begin{equation}\label{MP3}
\mbox{ the C-matrix of all contrasts for the plan } \tilde{C} ({\cal P}_0) = 6I_3, \end{equation}
so that condition (b) of the same theorem is also satisfied.
Hence the result follows from the same theorem. $\Box$

\vspace{.5em}

Next we generate an infinite series of POTBs from ${\cal P}_0$.

\begin{theo}\label{3^{3N}} (a) If an OA(N,m-1,3,2) exists, then
 there exists a {\bf connected} POTB ${\cal P}_m$ for a $3^{3m}$
experiment in $N$ blocks. Among the blocks, $2N$ are of size four,
 while the remaining $N$ are  of size $2$ each.

(b)  ${\cal P}_m$ is universally optimal among all plans for the same
experiment on the same set up.
\end{theo}

{\bf Proof :} (a) By assumption $Q(N,m,s)$ exists (see Notation \ref {ExtOA}). Let ${\cal P}_m = Q(N,m,s) \Diamond {\cal P}_0$,
where ${\cal P}_0$ is as in Table 3.2.3.
Now,  Theorem \ref {PlanThruOA} together with Theorem \ref {3^3} imply that the plan ${\cal P}_m$ is a POTB. That the block sizes are as in the statement follows from the construction.

  (b) Theorem \ref {PlanThruOA} together with \ref {3^3} imply that the information matrix of ${\cal P}_h$ is $6NI_{3m}$. The rest are
 as in the proof of Theorem \ref {3^3}. $\Box$

{\bf A series of plans for asymmetrical experiments.}

Let $s = 2t + 1$ be a prime power. As a special case of Lemma 2.6 of Morgan and Uddin (1996) we get a main effect plan for an $s^t$ experiment on $2s$ blocks of size $t = 1$ each. We shall add a factor with $s + 1$ levels to get a  main effect plan for an asymmetrical experiment. We follow  Notation \ref {sPrime} (a). Further, $F^+$ will denote the set $F \cup \{\infty\}$, where $F$ is as in Notation \ref {sPrime}. The following rule will define addition in $F^+$.
 \begin{equation}\label{addInfty} \infty + u = u + \infty = \infty , \; u \in F.\end{equation}

We define an $s \times s$ matrix $M$ with  rows and columns indexed by $F$.
 \begin{equation}\label{IncMatM} M = ((m (p,q)))_{p,q \in F}, \mbox{ where } m(p,q) = \left\{ \begin{array}{ll}
1 & \mbox{ if }  q - p  \in  C_0, \\
 0 & \mbox{ otherwise } \end{array} \right.\end{equation}

\begin{theo}\label{level(s+1)} Suppose  $s$ is a prime power. Let $t= (s-1)/2$.
Then we have the following.

(a) There exists a POTB ${\cal P}^*$ for an $s^t(s+1)$ experiment on $2s$ blocks of size $t+1$ each.

(b) ${\cal P}^*$ is  universally optimal for the inference on each $s$-level factor.

(c)  If $s \equiv 3 \pmod 4$, then
${\cal P}^*$ is  universally optimal for the inference on the $s+1$-level factor too.

\end{theo}

{\bf Proof :}  Take $\alpha \in C_1$. Consider a pair of $(t+1) \times (t+1)$ arrays $E^l, l = 0,1$  with rows and columns indexed by $C_0 \cup\{\infty\}$. The entries of $E^0$ and $E^1$ are as follows. For $x,y \in C_0,\; l = 0,1$,
\begin{eqnarray*} e^l_{x,y} = \alpha^{l}xy, \;   &  e^l (i,\infty) = 0,  \\
e^l_{\infty,y} = \alpha^{l+1} y, \;  & e^0(\infty,\infty) = 0 \mbox{ and }  e^1(\infty,\infty) = \infty.\end{eqnarray*}
Consider an initial plan ${\cal P}$  consisting of two initial blocks, the runs in which are the columns of $E^0$ and $E^1$ respectively.
The factors of ${\cal P}$ are indexed in accordance with the rows of $E_0,E_1$; that is, ${\cal F} = C_0 \cup\{\infty\}$. Let  ${\cal P}^* = {\cal P} \oplus F$.

 Clearly, ${\cal P}^*$ is a plan for an $s^t$(s+1) experiment on $2s$ blocks of size $t+1$ each. We can verify the following properties of the incidence matrices of  ${\cal P}^*$. For $x \in C_0$, we have
 \begin{eqnarray*} N_{xy} = I_s + J_s,\; N_{x\infty} = J_{s \times s+1} \\
 \mbox{ and } L_x = \left [\begin{array}{cc} M + I & J - M \end{array} \right ]. \end{eqnarray*}
 Here $ M $   is as in (\ref {IncMatM}).
 It follows that for $x \neq y, x,y \in C_0$, $L_x L'_y = (t+1)N_{xy}$, so that  $x$ and $y$ satisfies (\ref {Orth-bl}) and hence $x \perp y$.
  Again, $N_{x,\infty}$ satisfies (\ref {PFC}), so that $\infty$   is orthogonal to $x$ in the usual sense, $x \in C_0$. So, (a) is proved.

Next, we see that for $x \in C_0$, $L_x$ is the incidence matrix of a BIBD $d_1$ with parameters
$ (v = s, b = 2s, r = s+1, k = t + 1,  \lambda = t+1)$. Thus,  the conditions of Theorem \ref {univOpt} (a) is satisfied for each $x \in C_0$ and hence  (b) is proved by the same theorem.

 Finally, we see that if $s \equiv 3 \pmod 4$, then  $L_\infty$ is the incidence matrix of a BIBD $d_2$ with parameters
$ (v = s+1, b = 2s, r = s, k = t + 1,  \lambda = t)$. Therefore,  the conditions of Theorem \ref {univOpt} (a) are  satisfied for the factor
$\infty$ and hence  (c) is proved by the same theorem. $\Box$

\section{Reference} \begin{enumerate}
\item Addelman, S. (1962). Orthogonal main effect plans for asymmetrical
factorial experiments, Technometrics, vol. 4, p: 21-46.


\item Bagchi B and Bagchi S (2001). Optimality of partial geometric designs. Ann. Stat., vol. 29, p : 577-594.


\item  Bagchi, S. (2010). Main effect plans orthogonal through the
block factor. Technometrics, vol. 52, p : 243-249.

\item Bagchi, S. (2019). Inter-class orthogonal main effect plans for
asymmetrical experiments,  Sankhya B, vol. 81, p: 93–122.


\item Bose, M. and  Bagchi, S. (2007). Optimal main effect plans in blocks
of small size, Jour. Stat. Prob. Let., vol.  77, p : 142-147.

\item Chen, X.P., JG Lin, J.G., Yang, J.F. and Wang, H.X (2015). Construction of main effects plans orthogonal through the block factor,
Statist. Probab. Lett. , vol. 106,  p :  58-64

\item Das, A. and Dey, A. (2004). Optimal main effect plans with nonorthogonal
blocks. Sankhya, vol. 66, p : 378-384.




\item Hedayat, A.S., Sloan, N.J.A. and Stufken, J. (1999). Orthogonal
arrays, Theory and Applications, Springer Series in Statistics.

\item  Huang, L., Wu, C.F.J. and Yen, C.H. (2002). The idle
column method : Design construction, properties and comparisons,
Technometrics, vol. 44, p : 347-368.



\item Jacroux, Mike (2011). On the D-optimality of orthogonal and nonorthogonal blocked main effects plans. Statist. Probab. Lett. Vol. 81 ,  p:  116-120.
\item Jacroux, Mike (2011). On the D-optimality of nonorthogonal blocked main effects plans. Sankhya B, vol.  73,  p:   62-69.
\item Jacroux, Mike (2013). A note on the optimality of 2-level main effects plans in blocks of odd size. Statist. Probab. Lett. , vol. 83,  p:   1163-1166.
\item Jacroux, Mike, Kealy-Dichone, Bonni (2014). On the E-optimality of blocked main effects plans when $n  \equiv 3 \pmod 4$. Statist. Probab. Lett., vol.  87 ,  p:  143-148.

  \item Jacroux, Mike, Kealy-Dichone, Bonni (2015). On the E-optimality of blocked main effects plans when  $n \equiv 2 \pmod 4$. Sankhya B vol 77 ,  p:   165-174.
\item Jacroux, Mike, Jacroux, Tom  (2016).  On the E-optimality of blocked main effects plans when  $n \equiv 1 \pmod 4$. Comm. Statist. Theory Methods, vol. 45 ,  p:   5584-5589.

 \item Jacroux, Mike, Kealy-Dichone, Bonni  (2017). On the E-optimality of blocked main effects plans in blocks of different sizes. Comm. Statist. Theory Methods, vol.  46,  p:    2132-2138

\item Morgan, J.P. and Uddin, N. (1996). Optimal blocked main effect
plans with nested rows and columns and related designs. Ann. Stat.
vol. 24,   p:  1185-1208.

\item Mukerjee, R., Dey, A. and Chatterjee, K. (2002). Optimal main effect plans with non-orthogonal blocking. Biometrika, 89,  p:  225-229.

\item Rao, C.R. (1946). On Hypercubes of strength d and a system of confounding in factorial experiments Bull. Cal. Math. Soc., 38, p: 67.
1946

\item SahaRay, R., Dutta, G. (2016). On the Optimality of Blocked Main Effects Plans. International Scholarly and Scientific Research and Innovation, 10, p : 583-586.

experiments in blocks of size two. Sankhya, vol. 66, p : 327

\item Yanai, H., Takeuchi, K. and  Takane, K. (2011). Projection Matrices, Generalized Inverse Matrices, and Singular Value Decomposition, Statistics for Social and Behavioral Sciences. Springer.

\end{enumerate}

* Foot note : Both the authors are retired from Indian Statistical Institute, Bangalore Center

 \end{document}